\documentclass{article}
\usepackage{amsmath}

\input{tcilatex}
\begin{document}

\title{{\LARGE Pandiagonal and Knut Vik Sudoku Squares\thanks{%
Reprinted (with minor revisions) from: Mathematics Today 49, 86-87 (2013).
Appendix added.}}}
\author{Ronald P. Nordgren\thinspace \thanks{%
email: nordgren@rice.edu \quad Address: 3989 Pebble Beach Dr, Longmont, CO
80503-8358} \\
Brown School of Engineering, Rice University\\
}
\date{}
\maketitle

\noindent In her informative paper on the construction of pandiagonal magic
squares, Dame Kathleen Ollerenshaw \cite{OLL} conjectured that there are no
Knut Vik (or Latin pandiagonal) sudoku squares of order 9. In a subsequent
letter, Boyer \cite{BOY} confirmed that there are no such squares of order $%
2k$ and $3k$ since the nonexistence of Knut Vik squares of such orders was
proved by Hedayat \cite{HED}. However, this fact does not preclude the
existence of pandiagonal sudoku squares. Here we present systematic methods
of constructing pandiagonal sudoku squares of order $k^{2}$ and Knut Vik
sudoku squares of order $k^{2}$ not divisible by $2$ or $3.\smallskip $

\begin{center}
{\large Definitions}
\end{center}

\noindent For our purposes, an order-$n$ \emph{Latin }square matrix has the
integer elements $0,1,\ldots ,n-1$ on all rows and all columns which thus
sum to the \emph{index}%
\begin{equation}
m=n\left( n-1\right) /2\,.  \label{m}
\end{equation}%
Thus, a Latin square also is \emph{semi-magic}. In a \emph{pandiagonal}
square all $2n$ broken diagonals (with wraparound) in both directions sum to 
$m.$ A \emph{Knut Vik} square is a Latin square that has the integers $%
0,1,\ldots ,n-1$ on all broken diagonals; thus, it also is pandiagonal. A 
\emph{sudoku} square is a Latin square of order $n=k^{2}$ whose $n$ main,
order-$k$ subsquares also contain the integers $0,1,\ldots ,n-1$. We define
a \emph{super-sudoku }square\emph{\ }as a sudoku square with the additional
property that all $n^{2}$ of its order-$k$ subsquares, including broken
ones, sum to $m.\smallskip $

\begin{center}
{\large Pandiagonal Sudoku Squares}
\end{center}

\noindent By the foregoing definitions, the nonexistence of a Knut Vik
sudoku square does not preclude the existence of a pandiagonal sudoku square
of the same order. For example, here is a pandiagonal super-sudoku square
matrix of order 9:

\begin{equation}
S_{9}=\left[ 
\begin{array}{ccccccccc}
0 & 1 & 2 & 3 & 4 & 5 & 6 & 7 & 8 \\ 
3 & 4 & 5 & 6 & 7 & 8 & 0 & 1 & 2 \\ 
6 & 7 & 8 & 0 & 1 & 2 & 3 & 4 & 5 \\ 
1 & 2 & 0 & 4 & 5 & 3 & 7 & 8 & 6 \\ 
4 & 5 & 3 & 7 & 8 & 6 & 1 & 2 & 0 \\ 
7 & 8 & 6 & 1 & 2 & 0 & 4 & 5 & 3 \\ 
2 & 0 & 1 & 5 & 3 & 4 & 8 & 6 & 7 \\ 
5 & 3 & 4 & 8 & 6 & 7 & 2 & 0 & 1 \\ 
8 & 6 & 7 & 2 & 0 & 1 & 5 & 3 & 4%
\end{array}%
\right]  \label{Su9}
\end{equation}%
which is formed by a simple row-wise permutation scheme. This square is not
Knut Vik (none of order-9 exist). It can be decomposed into two auxiliary
matrices as%
\begin{equation}
S_{9}=3\hat{S}_{9}+\hat{S}_{9}^{T},  \label{S9dec}
\end{equation}%
where 
\begin{equation}
\hat{S}_{9}=\left[ 
\begin{array}{ccccccccc}
0 & 0 & 0 & 1 & 1 & 1 & 2 & 2 & 2 \\ 
1 & 1 & 1 & 2 & 2 & 2 & 0 & 0 & 0 \\ 
2 & 2 & 2 & 0 & 0 & 0 & 1 & 1 & 1 \\ 
0 & 0 & 0 & 1 & 1 & 1 & 2 & 2 & 2 \\ 
1 & 1 & 1 & 2 & 2 & 2 & 0 & 0 & 0 \\ 
2 & 2 & 2 & 0 & 0 & 0 & 1 & 1 & 1 \\ 
0 & 0 & 0 & 1 & 1 & 1 & 2 & 2 & 2 \\ 
1 & 1 & 1 & 2 & 2 & 2 & 0 & 0 & 0 \\ 
2 & 2 & 2 & 0 & 0 & 0 & 1 & 1 & 1%
\end{array}%
\right] ,\quad \hat{S}_{9}^{T}=\left[ 
\begin{array}{ccccccccc}
0 & 1 & 2 & 0 & 1 & 2 & 0 & 1 & 2 \\ 
0 & 1 & 2 & 0 & 1 & 2 & 0 & 1 & 2 \\ 
0 & 1 & 2 & 0 & 1 & 2 & 0 & 1 & 2 \\ 
1 & 2 & 0 & 1 & 2 & 0 & 1 & 2 & 0 \\ 
1 & 2 & 0 & 1 & 2 & 0 & 1 & 2 & 0 \\ 
1 & 2 & 0 & 1 & 2 & 0 & 1 & 2 & 0 \\ 
2 & 0 & 1 & 2 & 0 & 1 & 2 & 0 & 1 \\ 
2 & 0 & 1 & 2 & 0 & 1 & 2 & 0 & 1 \\ 
2 & 0 & 1 & 2 & 0 & 1 & 2 & 0 & 1%
\end{array}%
\right] .  \label{S9aux}
\end{equation}%
The form of $\hat{S}_{9}$ shows that $S_{9}$ is a pandiagonal super-sudoku
square.

The permutation scheme used to construct the pandiagonal super-sudoku square
(\ref{Su9}) can be generalized to construct pandiagonal super-sudoku squares
of order $n=k^{2}$. Also, the simple form of $\hat{S}_{9}$ can be
generalized to $\hat{S}_{n}$ $\left( n=k^{2}\right) $ with elements $\hat{S}%
_{ij}^{\left( n\right) }$ given by%
\begin{equation}
\hat{S}_{ij}^{\left( n\right) }=\left( i-1+\left\lfloor \frac{j-1}{k}%
\right\rfloor \right) \func{mod}k,\quad i,j=1,2,\ldots ,n,\quad n=k^{2}.
\label{Snij}
\end{equation}%
Then%
\begin{equation}
S_{n}=k\hat{S}_{n}+\hat{S}_{n}^{T},  \label{Sn}
\end{equation}%
is a pandiagonal super-sudoku square of the same permutation form as $S_{9}.$
The pandiagonal and super-sudoku properties of $\hat{S}_{n}$ can be verified
using formulas given by Nordgren [4,5].

In addition, we observe that a natural pandiagonal magic square $M_{n}$ with
elements $0,1,\cdots ,n^{2}-1$ can be constructed from a pandiagonal
super-sudoku square $S_{n}$ of odd order $n$ by the auxiliary square formula 
\begin{equation}
M_{n}=nS_{n}+S_{n}R,  \label{Mn}
\end{equation}%
where $R$ is the reflection matrix with $1$ on its cross diagonal (top-right
to bottom-left) and $0$ for all other elements. Here, $S_{n}$ and its
reflection about its vertical centerline, $S_{n}R,$ are seen to be
orthogonal (as defined in [1,4]). For example, for $n=9$, from (\ref{Su9})
and (\ref{Mn}), we obtain 
\begin{equation}
M_{9}=\left[ 
\begin{array}{ccccccccc}
8 & 16 & 24 & 32 & 40 & 48 & 56 & 64 & 72 \\ 
29 & 37 & 45 & 62 & 70 & 78 & 5 & 13 & 21 \\ 
59 & 67 & 75 & 2 & 10 & 18 & 35 & 43 & 51 \\ 
15 & 26 & 7 & 39 & 50 & 31 & 63 & 74 & 55 \\ 
36 & 47 & 28 & 69 & 80 & 61 & 12 & 23 & 4 \\ 
66 & 77 & 58 & 9 & 20 & 1 & 42 & 53 & 34 \\ 
25 & 6 & 17 & 49 & 30 & 41 & 73 & 54 & 65 \\ 
46 & 27 & 38 & 79 & 60 & 71 & 22 & 3 & 14 \\ 
76 & 57 & 68 & 19 & 0 & 11 & 52 & 33 & 44%
\end{array}%
\right]  \label{M99}
\end{equation}%
in which all 3 by 3 subsquares, including broken ones, also sum to the index 
$360.\smallskip $

\begin{center}
{\large Knut Vik Sudoku Squares}
\end{center}

\noindent Boyer \cite{BOY} gives a Knut Vik sudoku square of order 25, the
lowest order for which such a square exists. However, he does not indicate
the method of constructing it nor how to construct such squares of higher
order. On decomposing Boyer's square, we find a general method of
constructing a Knut Vik sudoku square $V_{n}$ from auxiliary squares using%
\begin{equation}
V_{n}=kV_{n}^{\prime }+V_{n}^{\prime \prime }  \label{MauxG}
\end{equation}%
for order $n=k^{2}$ not divisible by $2$ or $3$. The first auxiliary matrix $%
V_{n}^{\prime }$ has the same submatrix $A_{k}$ for its $n$ sudoku
subsquares, where $A_{k}$ is the order-$k$ submatrix with integers $%
0,1,\ldots ,k-1$ on its main diagonal which is then replicated via a chess
knight's move of right two/down one (with wraparound), e.g., for $n=25,$ $%
k=5 $%
\begin{equation}
V_{25}^{\prime }=\left[ 
\begin{array}{lllll}
A_{5} & A_{5} & A_{5} & A_{5} & A_{5} \\ 
A_{5} & A_{5} & A_{5} & A_{5} & A_{5} \\ 
A_{5} & A_{5} & A_{5} & A_{5} & A_{5} \\ 
A_{5} & A_{5} & A_{5} & A_{5} & A_{5} \\ 
A_{5} & A_{5} & A_{5} & A_{5} & A_{5}%
\end{array}%
\right] ,\quad A_{5}=\left[ 
\begin{array}{ccccc}
0 & 4 & 3 & 2 & 1 \\ 
2 & 1 & 0 & 4 & 3 \\ 
4 & 3 & 2 & 1 & 0 \\ 
1 & 0 & 4 & 3 & 2 \\ 
3 & 2 & 1 & 0 & 4%
\end{array}%
\right] .  \label{A}
\end{equation}%
The form of the second auxiliary matrix $V_{n}^{\prime \prime }$ is
illustrated for $n=25,$ $k=5$ by%
\begin{equation}
V_{5}^{\prime \prime }=\left[ 
\begin{array}{ccccc}
B & E & C & F & D \\ 
C & F & D & B & E \\ 
D & B & E & C & F \\ 
E & C & F & D & B \\ 
F & D & B & E & C%
\end{array}%
\right] ,  \label{Mt5}
\end{equation}%
where the submatrices $B,C,D,E,F$ replicate in a chess knight's move of left
two/down one (with wraparound) and are given by%
\begin{equation}
B=\left[ 
\begin{array}{ccccc}
0 & 3 & 1 & 4 & 2 \\ 
2 & 0 & 3 & 1 & 4 \\ 
4 & 2 & 0 & 3 & 1 \\ 
1 & 4 & 2 & 0 & 3 \\ 
3 & 1 & 4 & 2 & 0%
\end{array}%
\right] ,\quad C=\left[ 
\begin{array}{ccccc}
2 & 0 & 3 & 1 & 4 \\ 
4 & 2 & 0 & 3 & 1 \\ 
1 & 4 & 2 & 0 & 3 \\ 
3 & 1 & 4 & 2 & 0 \\ 
0 & 3 & 1 & 4 & 2%
\end{array}%
\right] ,\quad \cdots ,  \label{M5sub}
\end{equation}%
Here, $B$ is formed by starting its first column with the first column of $%
A_{5}$ and replicating via a chess bishop's move of down one/right one (with
wraparound). Next, $C$ is formed by shifting the rows of $B$ up one row (top
row to bottom) and similarly for $C\Rightarrow D\Rightarrow E\Rightarrow F$
which form the first column of $V_{5}^{\prime \prime }$. The second column
of $V_{5}^{\prime \prime }$ starts with submatrix $E$ which has the element $%
1$ in its upper-left corner. The third column of $V_{5}^{\prime \prime }$
starts with submatrix $C$ which has the element $2$ in its upper-left corner
and so on. The submatrices in each column progress downward in the same
order as in the first column (with wraparound).

Higher-order Knut Vik sudoku squares can be constructed in a similar manner
with the elements of $V_{n}$ given by%
\begin{gather}
V_{n}\left( i,j\right) =k\left( \left( 2i-j-1\right) \func{mod}k\right)
+\left( 2i+2\left\lfloor \frac{i-1}{k}\right\rfloor -2j+\left\lfloor \frac{%
j-1}{k}\right\rfloor \right) \func{mod}k,  \label{Kij} \\
i,j=1,2,\ldots ,n,\quad n=k^{2},\quad n\func{mod}2\neq 0,\quad n\func{mod}%
3\neq 0\,.  \notag
\end{gather}%
The order-$k$ sudoku squares again form a Knight's move pattern of left
two/down one (with wraparound). A detailed analysis shows that $V_{n}$ from (%
\ref{Kij}) is a Knut Vik sudoku square. As in the construction of
pandiagonal magic squares given in \cite{OLL} and \cite{NOR1}, the knight's
move patterns of $A_{k}$ and $V_{n}^{\prime \prime }$, as embodied in (\ref%
{Kij}), do not work when the order $n=k^{2}$ is a multiple of $2$ or $3,$
i.e., when a Knut Vik square does not exist.

An order-25, Knut Vik sudoku square of the form given by Boyer \cite{BOY}
follows from (\ref{Kij}). Also, a natural pandiagonal magic square can be
constructed from a Knut Vik sudoku square using (\ref{Mn}). Example
pandiagonal and Knut Vik sudoku squares constructed by the methods presented
here are given in the Appendix.

\begin{center}
$\smallskip $\newpage

{\Large Appendix: Examples}$\smallskip $
\end{center}

\noindent Here we present more examples of pandiagonal and Knut Vik sudoku
squares constructed using Excel\copyright\ and Maple\copyright\ in the
Scientific WorkPlace\copyright\ program. The special properties of these
squares are verified using defining formulas from \cite{NOR1} or \cite{NOR2}
which are reviewed next.$\smallskip $

\begin{center}
{\large Definitions}
\end{center}

\noindent First we define two permutation matrices that are used in what
follows. Let $R$ denote the order-$n$ \emph{reflection matrix} with $1$ on
its cross diagonal (top-right to bottom-left) and $0$ for all other
elements. For a given order-$n$ matrix $M,$ the operation $RM$ reflects the
elements of $M$ about its horizontal centerline and $MR$ reflects the
elements of $M$ about its vertical centerline.

Let $K$ denote the order-$n$ \emph{shifter matrix} that has all elements $0$
except $K\left( 1,n\right) =1$ (upper right corner) and $K\left(
i,i-1\right) =1\,,\ i=2,3,\cdots ,n$ (diagonal below the main diagonal). It
can be defined by%
\begin{equation}
K(i,j)=\delta \left( \left( i-j+n-1\right) \func{mod}n\right) ,\quad
i,j=1,2,\ldots ,n,  \label{Kijn}
\end{equation}%
where $\delta \left( k\right) $ is the Dirac delta for integers, defined as%
\begin{equation}
\delta \left( k\right) =\left\{ 
\begin{tabular}{c}
$1$ if $k=0$ \\ 
$0$ if $k\neq 0$%
\end{tabular}%
\right.
\end{equation}%
and in Maple\copyright\ as the function%
\begin{equation}
\delta \left( k\right) =\func{Heaviside}\left( k+\frac{1}{2}\right) -\func{%
Heaviside}\left( k-\frac{1}{2}\right) .
\end{equation}%
The operation $KM$ shifts rows of $M$ down one (and bottom row to top) while 
$MK$ shifts columns of $M$ one to the left (and first column to last). Power
operations $K^{i}M$ and $MK^{i}$ give rise to repeated shifts. The following
identities can be easily verified:%
\begin{gather}
K^{n}=I,\quad K^{-i}=K^{n-i},\quad i=1,2,\ldots ,n,  \notag \\
\sum_{i=1}^{n}K^{i}=\sum_{i=1}^{n}K^{-i}=U.  \label{KiU}
\end{gather}

An order-$n$ matrix $M$ is \emph{semi-magic\ }if the sums of all its rows
and all its columns equal the same \emph{index} $m$, i.e., if%
\begin{equation}
Mu=\left( u^{T}M\right) ^{T}=mu\quad \text{or}\quad MU=UM=mU,  \label{UM}
\end{equation}%
where $u$ is an order-$n$ column vector with all elements $1$ and $U$ is an
order-$n$ matrix with all elements $1\,.$ The matrix $M$ is \emph{magic\ }if
in addition to (\ref{UM}) its main diagonal and cross diagonal also sum to $%
m,$ i.e., if%
\begin{equation}
\limfunc{tr}\left[ M\right] =\limfunc{tr}\left[ RM\right] =m.  \label{TrM}
\end{equation}%
By these definitions, a magic square also is semi-magic and a semi-magic
square may or may not be magic. Subscripts are used to denote special
classes of $M.$ The matrix $M_{N}$ is \emph{natural} if its elements are
integers in the numerical sequence $0,1,\ldots ,n^{2}-1.$ The natural
property of $M_{N}$ can be verified by sorting its elements into numerical
order or from checking that%
\begin{equation}
M_{N}-K^{i}M_{N}K^{j}=\Omega ,\quad i,j=1,2,\ldots ,n,\quad \left( i=j\neq
n\right) ,  \label{Nat}
\end{equation}%
where $\Omega $ is any order-$n$ matrix that has no zero elements. The index
of a natural magic or semi-magic square is%
\begin{equation}
m_{N}=n(n^{2}-1)/2\,.  \label{IndexN}
\end{equation}

An order-$n$ square matrix $M_{P}$ is \emph{pandiagonal} if all its broken
diagonals (of $n$ elements) in both directions sum to the index $m,$ i.e., if%
\begin{equation}
\limfunc{tr}\left[ K^{i}M_{P}\right] =\limfunc{tr}\left[ K^{i}RM_{P}\right]
=m,\quad i=1,2,\ldots ,n.
\end{equation}%
Also, the pandiagonal property of $M_{P}$ can be verified from the identities%
\begin{equation}
\sum_{i=1}^{n}K^{i}M_{P}K^{i}=mU\quad \text{and}\quad
\sum_{i=1}^{n}K^{i}M_{P}K^{-i}=mU.  \label{KMK}
\end{equation}
A square that is pandiagonal and magic is called \emph{panmagic}.

An order-$n$ \emph{Latin }square matrix $M_{L}$ has the integer elements $%
0,1,\ldots ,n-1$ on all rows and all columns which thus sum to the index%
\begin{equation}
m=n\left( n-1\right) /2\,.  \label{Index}
\end{equation}%
Thus, a Latin square also is semi-magic. The Latin square property of $M_{L}$
can be verified by sorting its rows and columns into numerical order or from
checking that%
\begin{equation}
M_{L}-K^{i}M_{L}=\Omega ,\quad M_{L}-M_{L}K^{i}=\Omega ,\quad i=1,2,\ldots
,n-1\,.  \label{Lat}
\end{equation}

An order-$n$ \emph{Knut Vik} square $V$ is a Latin square that has the
integers $0,1,\ldots ,n-1$ on all broken diagonals; thus, it also is
pandiagonal. The Knut Vik diagonal property of $V$ can be verified by
sorting its diagonals into numerical order or by checking that 
\begin{equation}
V-K^{i}VK^{i}=\Omega ,\quad V-K^{i}VK^{-i}=\Omega ,\quad i=1,2,\ldots ,n-1\,.
\label{KV}
\end{equation}

A \emph{sudoku} square $S$ is a Latin square of order $n=k^{2}$ whose $n$
main, order-$k$ subsquares also contain the integers $0,1,\ldots ,n-1$. The
sudoku\emph{\ }property of a square can be verified by sorting the elements
of its sudoku subsquares or by checking a matrix equation formed by first
defining the order-$n$ permutation matrix $H$ in terms of submatrices $K_{k}$
and $O_{k}$ as%
\begin{equation}
H=\left[ 
\begin{array}{cccc}
K_{k} & O_{k} & \cdots & O_{k} \\ 
O_{k} & \ddots &  & O_{k} \\ 
\vdots &  & \ddots & \vdots \\ 
O_{k} & O_{k} & \cdots & K_{k}%
\end{array}%
\right] ,  \label{Hn}
\end{equation}%
where $K_{k}$ is an order-$k$ shifter matrix and $O_{k}$ is an order-$k$
matrix with all elements $0$. Then, $S$ is a Sudoku square if 
\begin{equation}
S-H^{i}SH^{j}=\Omega ,\quad i,j=1,2,\ldots ,k,\quad \left( i=j\neq k\right) .
\label{Sud}
\end{equation}

We define a \emph{super-sudoku }square\emph{\ }as a sudoku square with the
additional property that all $n^{2}$ of its order-$k$ subsquares, including
broken ones, sum to $m.$ The super-sudoku sum property of $S$ can be
verified from%
\begin{equation}
\left( I+K+K^{2}+\cdots +K^{k-1}\right) S\left( I+K+K^{2}+\cdots
+K^{k-1}\right) =mU.  \label{KSK}
\end{equation}

\begin{center}
{\large Pandiagonal Sudoku Squares}
\end{center}

\noindent For $n=9,$ by (\ref{Snij}) and (\ref{Sn}), we have\footnote{%
{\footnotesize Use the\textquotedblleft New Definition\textquotedblright\
command in Maple\copyright\ to define }$S_{9}(i,j)${\footnotesize \ and use
the \textquotedblleft Fill Matrix\textquotedblright\ command to form it.}}%
\begin{gather}
S_{9}\left( i,j\right) =3\left( \left( i-1+\left\lfloor \frac{j-1}{3}%
\right\rfloor \right) \func{mod}3\right)  \notag \\
+\left( \left( j-1+\left\lfloor \frac{i-1}{3}\right\rfloor \right) \func{mod}%
3\right) \equiv s_{9},\quad i,j=1,2,\ldots ,9,  \label{s9} \\
s_{9}=\left[ 
\begin{array}{ccccccccc}
0 & 1 & 2 & 3 & 4 & 5 & 6 & 7 & 8 \\ 
3 & 4 & 5 & 6 & 7 & 8 & 0 & 1 & 2 \\ 
6 & 7 & 8 & 0 & 1 & 2 & 3 & 4 & 5 \\ 
1 & 2 & 0 & 4 & 5 & 3 & 7 & 8 & 6 \\ 
4 & 5 & 3 & 7 & 8 & 6 & 1 & 2 & 0 \\ 
7 & 8 & 6 & 1 & 2 & 0 & 4 & 5 & 3 \\ 
2 & 0 & 1 & 5 & 3 & 4 & 8 & 6 & 7 \\ 
5 & 3 & 4 & 8 & 6 & 7 & 2 & 0 & 1 \\ 
8 & 6 & 7 & 2 & 0 & 1 & 5 & 3 & 4%
\end{array}%
\right]  \notag
\end{gather}%
and (\ref{Kijn}) gives%
\begin{align}
K_{9}(i,j)& =\delta \left( \left( i-j+9-1\right) \func{mod}9\right) \equiv
k_{9},\quad i,j=1,2,\ldots ,9, \\
k_{9}& =\left[ 
\begin{array}{ccccccccc}
0 & 0 & 0 & 0 & 0 & 0 & 0 & 0 & 1 \\ 
1 & 0 & 0 & 0 & 0 & 0 & 0 & 0 & 0 \\ 
0 & 1 & 0 & 0 & 0 & 0 & 0 & 0 & 0 \\ 
0 & 0 & 1 & 0 & 0 & 0 & 0 & 0 & 0 \\ 
0 & 0 & 0 & 1 & 0 & 0 & 0 & 0 & 0 \\ 
0 & 0 & 0 & 0 & 1 & 0 & 0 & 0 & 0 \\ 
0 & 0 & 0 & 0 & 0 & 1 & 0 & 0 & 0 \\ 
0 & 0 & 0 & 0 & 0 & 0 & 1 & 0 & 0 \\ 
0 & 0 & 0 & 0 & 0 & 0 & 0 & 1 & 0%
\end{array}%
\right] .  \notag
\end{align}%
$\allowbreak \allowbreak $The pandiagonal condition (\ref{KMK}) for $s_{9}$
can be evaluated using Maple\copyright\ as follows:%
\begin{equation}
\sum_{r=1}^{9}\left[ \left[ x^{r}\right] _{x=k_{9}}\left[ x\right] _{x=s_{9}}%
\right] \left[ x^{r}\right] _{x=k_{9}}=\sum_{r=1}^{9}\left[ \left[ x^{r}%
\right] _{x=k_{9}}\left[ x\right] _{x=s_{9}}\right] \left[ x^{9-r}\right]
_{x=k_{9}}=36U
\end{equation}%
which verifies that $s_{9}$ is pandiagonal. The sudoku condition (\ref{Sud})
can be verified with $h_{9}$ from (\ref{Hn}) as%
\begin{equation}
h_{9}=\left[ 
\begin{array}{ccccccccc}
0 & 0 & 1 & 0 & 0 & 0 & 0 & 0 & 0 \\ 
1 & 0 & 0 & 0 & 0 & 0 & 0 & 0 & 0 \\ 
0 & 1 & 0 & 0 & 0 & 0 & 0 & 0 & 0 \\ 
0 & 0 & 0 & 0 & 0 & 1 & 0 & 0 & 0 \\ 
0 & 0 & 0 & 1 & 0 & 0 & 0 & 0 & 0 \\ 
0 & 0 & 0 & 0 & 1 & 0 & 0 & 0 & 0 \\ 
0 & 0 & 0 & 0 & 0 & 0 & 0 & 0 & 1 \\ 
0 & 0 & 0 & 0 & 0 & 0 & 1 & 0 & 0 \\ 
0 & 0 & 0 & 0 & 0 & 0 & 0 & 1 & 0%
\end{array}%
\right] .
\end{equation}%
The verification of (\ref{Sud}) can be done in Maple\copyright\ by examining
each term in the sum%
\begin{equation}
\sum_{p=1}^{3}\sum_{r=1}^{3}\left( \left[ x\right] _{x=s_{9}}-\left[ \left[
x^{p}\right] _{x=h_{9}}\left[ x\right] _{x=s_{9}}\right] \left[ x^{r}\right]
_{x=h_{9}}\right) .
\end{equation}%
The super-sudoku sum equation (\ref{KSK}) can be evaluated using Maple%
\copyright\ as%
\begin{equation}
\left[ \left[ 1+x+x^{2}\right] _{x=k_{9}}\left[ x\right] _{x=s_{9}}\right] %
\left[ 1+x+x^{2}\right] _{x=k_{9}}=36U
\end{equation}%
which, together with the sudoku property, verifies that $s_{9}$ is
super-sudoku.

A natural panmagic square with elements $0,1,\ldots ,80$ can be constructed
from (\ref{Mn}) as%
\begin{align}
M_{9}\left( i,j\right) & =9S_{9}\left( i,j\right) +S_{9}\left( i,10-j\right)
\equiv m_{9},\quad i,j=1,2,\ldots ,9, \\
m_{9}& =\left[ 
\begin{array}{ccccccccc}
8 & 16 & 24 & 32 & 40 & 48 & 56 & 64 & 72 \\ 
29 & 37 & 45 & 62 & 70 & 78 & 5 & 13 & 21 \\ 
59 & 67 & 75 & 2 & 10 & 18 & 35 & 43 & 51 \\ 
15 & 26 & 7 & 39 & 50 & 31 & 63 & 74 & 55 \\ 
36 & 47 & 28 & 69 & 80 & 61 & 12 & 23 & 4 \\ 
66 & 77 & 58 & 9 & 20 & 1 & 42 & 53 & 34 \\ 
25 & 6 & 17 & 49 & 30 & 41 & 73 & 54 & 65 \\ 
46 & 27 & 38 & 79 & 60 & 71 & 22 & 3 & 14 \\ 
76 & 57 & 68 & 19 & 0 & 11 & 52 & 33 & 44%
\end{array}%
\right] .  \notag
\end{align}%
Checking properties of $M_{9}\left( i,j\right) $, we find that%
\begin{gather}
\sum_{i=1}^{9}M_{9}\left( i,j\right) =360,\quad j=1,\ldots ,9,  \notag \\
\sum_{j=1}^{9}M_{9}\left( i,j\right) =360,\quad i=1,\ldots ,9,  \notag \\
\sum_{r=1}^{9}\left[ \left[ x^{r}\right] _{x=k_{9}}\left[ x\right] _{x=m_{9}}%
\right] \left[ x^{r}\right] _{x=k_{9}}=360U, \\
\sum_{r=1}^{9}\left[ \left[ x^{r}\right] _{x=k_{9}}\left[ x\right] _{x=m_{9}}%
\right] \left[ x^{9-r}\right] _{x=k_{9}}=360U,  \notag \\
\left[ \left[ 1+x+x^{2}\right] _{x=k_{9}}\left[ x\right] _{x=m_{9}}\right] %
\left[ 1+x+x^{2}\right] _{x=k_{9}}=360U  \notag
\end{gather}%
which verifies that $M_{9}$ is a panmagic square whose 81 3 by 3 subsquares
all add to $m=360\,.$ The natural property can be verified from (\ref{Nat})
by examining each term in the sum%
\begin{equation}
\sum_{p=1}^{9}\sum_{r=1}^{9}\left( \left[ x\right] _{x=m_{9}}-\left[ \left[
x^{p}\right] _{x=h_{9}}\left[ x\right] _{x=m_{9}}\right] \left[ x^{r}\right]
_{x=h_{9}}\right) .
\end{equation}

For $n=16,$ by (\ref{Snij}) and (\ref{Sn}), we have%
\begin{equation}
T_{16}(i,j)=\left( i-1+\left\lfloor \frac{j-1}{4}\right\rfloor \right) \func{%
mod}4\equiv t_{16},\quad i,j=1,2,\ldots ,16,
\end{equation}%
\begin{equation*}
t_{16}=\left[ 
\begin{array}{cccccccccccccccc}
0 & 0 & 0 & 0 & 1 & 1 & 1 & 1 & 2 & 2 & 2 & 2 & 3 & 3 & 3 & 3 \\ 
1 & 1 & 1 & 1 & 2 & 2 & 2 & 2 & 3 & 3 & 3 & 3 & 0 & 0 & 0 & 0 \\ 
2 & 2 & 2 & 2 & 3 & 3 & 3 & 3 & 0 & 0 & 0 & 0 & 1 & 1 & 1 & 1 \\ 
3 & 3 & 3 & 3 & 0 & 0 & 0 & 0 & 1 & 1 & 1 & 1 & 2 & 2 & 2 & 2 \\ 
0 & 0 & 0 & 0 & 1 & 1 & 1 & 1 & 2 & 2 & 2 & 2 & 3 & 3 & 3 & 3 \\ 
1 & 1 & 1 & 1 & 2 & 2 & 2 & 2 & 3 & 3 & 3 & 3 & 0 & 0 & 0 & 0 \\ 
2 & 2 & 2 & 2 & 3 & 3 & 3 & 3 & 0 & 0 & 0 & 0 & 1 & 1 & 1 & 1 \\ 
3 & 3 & 3 & 3 & 0 & 0 & 0 & 0 & 1 & 1 & 1 & 1 & 2 & 2 & 2 & 2 \\ 
0 & 0 & 0 & 0 & 1 & 1 & 1 & 1 & 2 & 2 & 2 & 2 & 3 & 3 & 3 & 3 \\ 
1 & 1 & 1 & 1 & 2 & 2 & 2 & 2 & 3 & 3 & 3 & 3 & 0 & 0 & 0 & 0 \\ 
2 & 2 & 2 & 2 & 3 & 3 & 3 & 3 & 0 & 0 & 0 & 0 & 1 & 1 & 1 & 1 \\ 
3 & 3 & 3 & 3 & 0 & 0 & 0 & 0 & 1 & 1 & 1 & 1 & 2 & 2 & 2 & 2 \\ 
0 & 0 & 0 & 0 & 1 & 1 & 1 & 1 & 2 & 2 & 2 & 2 & 3 & 3 & 3 & 3 \\ 
1 & 1 & 1 & 1 & 2 & 2 & 2 & 2 & 3 & 3 & 3 & 3 & 0 & 0 & 0 & 0 \\ 
2 & 2 & 2 & 2 & 3 & 3 & 3 & 3 & 0 & 0 & 0 & 0 & 1 & 1 & 1 & 1 \\ 
3 & 3 & 3 & 3 & 0 & 0 & 0 & 0 & 1 & 1 & 1 & 1 & 2 & 2 & 2 & 2%
\end{array}%
\right] ,
\end{equation*}%
\begin{gather}
S_{16}\left( i,j\right) =4T_{16}(i,j)+T_{16}(j,i)=4t_{16}+t_{16}^{T}\equiv
s_{16}, \\
i,j=1,2,\ldots ,16,  \notag \\
s_{16}=\left[ \text{{\scriptsize \negthinspace }}%
\begin{array}{cccccccccccccccc}
\text{{\small 0}} & \text{{\small 1}} & \text{{\small 2}} & \text{{\small 3}}
& \text{{\small 4}} & \text{{\small 5}} & \text{{\small 6}} & \text{{\small 7%
}} & \text{{\small 8}} & \text{{\small 9}} & \text{{\small \negthinspace 10}}
& \text{{\small \negthinspace 11}} & \text{{\small \negthinspace 12}} & 
\text{{\small \negthinspace 13}} & \text{{\small \negthinspace 14}} & \text{%
{\small \negthinspace 15}} \\ 
\text{{\small 4}} & \text{{\small 5}} & \text{{\small 6}} & \text{{\small 7}}
& \text{{\small 8}} & \text{{\small 9}} & \text{{\small \negthinspace 10}} & 
\text{{\small \negthinspace 11}} & \text{{\small \negthinspace 12}} & \text{%
{\small \negthinspace 13}} & \text{{\small \negthinspace 14}} & \text{%
{\small \negthinspace 15}} & \text{{\small 0}} & \text{{\small 1}} & \text{%
{\small 2}} & \text{{\small 3}} \\ 
\text{{\small 8}} & \text{{\small 9}} & \text{{\small \negthinspace 10}} & 
\text{{\small \negthinspace 11}} & \text{{\small \negthinspace 12}} & \text{%
{\small \negthinspace 13}} & \text{{\small \negthinspace 14}} & \text{%
{\small \negthinspace 15}} & \text{{\small 0}} & \text{{\small 1}} & \text{%
{\small 2}} & \text{{\small 3}} & \text{{\small 4}} & \text{{\small 5}} & 
\text{{\small 6}} & \text{{\small 7}} \\ 
\text{{\small \negthinspace 12}} & \text{{\small \negthinspace 13}} & \text{%
{\small \negthinspace 14}} & \text{{\small \negthinspace 15}} & \text{%
{\small 0}} & \text{{\small 1}} & \text{{\small 2}} & \text{{\small 3}} & 
\text{{\small 4}} & \text{{\small 5}} & \text{{\small 6}} & \text{{\small 7}}
& \text{{\small 8}} & \text{{\small 9}} & \text{{\small \negthinspace 10}} & 
\text{{\small \negthinspace 11}} \\ 
\text{{\small 1}} & \text{{\small 2}} & \text{{\small 3}} & \text{{\small 0}}
& \text{{\small 5}} & \text{{\small 6}} & \text{{\small 7}} & \text{{\small 4%
}} & \text{{\small 9}} & \text{{\small \negthinspace 10}} & \text{{\small %
\negthinspace 11}} & \text{{\small 8}} & \text{{\small \negthinspace 13}} & 
\text{{\small \negthinspace 14}} & \text{{\small \negthinspace 15}} & \text{%
{\small \negthinspace 12}} \\ 
\text{{\small 5}} & \text{{\small 6}} & \text{{\small 7}} & \text{{\small 4}}
& \text{{\small 9}} & \text{{\small \negthinspace 10}} & \text{{\small %
\negthinspace 11}} & \text{{\small 8}} & \text{{\small \negthinspace 13}} & 
\text{{\small \negthinspace 14}} & \text{{\small \negthinspace 15}} & \text{%
{\small \negthinspace 12}} & \text{{\small 1}} & \text{{\small 2}} & \text{%
{\small 3}} & \text{{\small 0}} \\ 
\text{{\small 9}} & \text{{\small \negthinspace 10}} & \text{{\small %
\negthinspace 11}} & \text{{\small 8}} & \text{{\small \negthinspace 13}} & 
\text{{\small \negthinspace 14}} & \text{{\small \negthinspace 15}} & \text{%
{\small \negthinspace 12}} & \text{{\small 1}} & \text{{\small 2}} & \text{%
{\small 3}} & \text{{\small 0}} & \text{{\small 5}} & \text{{\small 6}} & 
\text{{\small 7}} & \text{{\small 4}} \\ 
\text{{\small \negthinspace 13}} & \text{{\small \negthinspace 14}} & \text{%
{\small \negthinspace 15}} & \text{{\small \negthinspace 12}} & \text{%
{\small 1}} & \text{{\small 2}} & \text{{\small 3}} & \text{{\small 0}} & 
\text{{\small 5}} & \text{{\small 6}} & \text{{\small 7}} & \text{{\small 4}}
& \text{{\small 9}} & \text{{\small \negthinspace 10}} & \text{{\small %
\negthinspace 11}} & \text{{\small 8}} \\ 
\text{{\small 2}} & \text{{\small 3}} & \text{{\small 0}} & \text{{\small 1}}
& \text{{\small 6}} & \text{{\small 7}} & \text{{\small 4}} & \text{{\small 5%
}} & \text{{\small \negthinspace 10}} & \text{{\small \negthinspace 11}} & 
\text{{\small 8}} & \text{{\small 9}} & \text{{\small \negthinspace 14}} & 
\text{{\small \negthinspace 15}} & \text{{\small \negthinspace 12}} & \text{%
{\small \negthinspace 13}} \\ 
\text{{\small 6}} & \text{{\small 7}} & \text{{\small 4}} & \text{{\small 5}}
& \text{{\small \negthinspace 10}} & \text{{\small \negthinspace 11}} & 
\text{{\small 8}} & \text{{\small 9}} & \text{{\small \negthinspace 14}} & 
\text{{\small \negthinspace 15}} & \text{{\small \negthinspace 12}} & \text{%
{\small \negthinspace 13}} & \text{{\small 2}} & \text{{\small 3}} & \text{%
{\small 0}} & \text{{\small 1}} \\ 
\text{{\small \negthinspace 10}} & \text{{\small \negthinspace 11}} & \text{%
{\small 8}} & \text{{\small 9}} & \text{{\small \negthinspace 14}} & \text{%
{\small \negthinspace 15}} & \text{{\small \negthinspace 12}} & \text{%
{\small \negthinspace 13}} & \text{{\small 2}} & \text{{\small 3}} & \text{%
{\small 0}} & \text{{\small 1}} & \text{{\small 6}} & \text{{\small 7}} & 
\text{{\small 4}} & \text{{\small 5}} \\ 
\text{{\small \negthinspace 14}} & \text{{\small \negthinspace 15}} & \text{%
{\small \negthinspace 12}} & \text{{\small \negthinspace 13}} & \text{%
{\small 2}} & \text{{\small 3}} & \text{{\small 0}} & \text{{\small 1}} & 
\text{{\small 6}} & \text{{\small 7}} & \text{{\small 4}} & \text{{\small 5}}
& \text{{\small \negthinspace 10}} & \text{{\small \negthinspace 11}} & 
\text{{\small 8}} & \text{{\small 9}} \\ 
\text{{\small 3}} & \text{{\small 0}} & \text{{\small 1}} & \text{{\small 2}}
& \text{{\small 7}} & \text{{\small 4}} & \text{{\small 5}} & \text{{\small 6%
}} & \text{{\small \negthinspace 11}} & \text{{\small 8}} & \text{{\small 9}}
& \text{{\small \negthinspace 10}} & \text{{\small \negthinspace 15}} & 
\text{{\small \negthinspace 12}} & \text{{\small \negthinspace 13}} & \text{%
{\small \negthinspace 14}} \\ 
\text{{\small 7}} & \text{{\small 4}} & \text{{\small 5}} & \text{{\small 6}}
& \text{{\small \negthinspace 11}} & \text{{\small 8}} & \text{{\small 9}} & 
\text{{\small \negthinspace 10}} & \text{{\small \negthinspace 15}} & \text{%
{\small \negthinspace 12}} & \text{{\small \negthinspace 13}} & \text{%
{\small \negthinspace 14}} & \text{{\small 3}} & \text{{\small 0}} & \text{%
{\small 1}} & \text{{\small 2}} \\ 
\text{{\small \negthinspace 11}} & \text{{\small 8}} & \text{{\small 9}} & 
\text{{\small \negthinspace 10}} & \text{{\small \negthinspace 15}} & \text{%
{\small \negthinspace 12}} & \text{{\small \negthinspace 13}} & \text{%
{\small \negthinspace 14}} & \text{{\small 3}} & \text{{\small 0}} & \text{%
{\small 1}} & \text{{\small 2}} & \text{{\small 7}} & \text{{\small 4}} & 
\text{{\small 5}} & \text{{\small 6}} \\ 
\text{{\small \negthinspace 15}} & \text{{\small \negthinspace 12}} & \text{%
{\small \negthinspace 13}} & \text{{\small \negthinspace 14}} & \text{%
{\small 3}} & \text{{\small 0}} & \text{{\small 1}} & \text{{\small 2}} & 
\text{{\small 7}} & \text{{\small 4}} & \text{{\small 5}} & \text{{\small 6}}
& \text{{\small \negthinspace 11}} & \text{{\small 8}} & \text{{\small 9}} & 
\text{{\small \negthinspace 10}}%
\end{array}%
\text{{\scriptsize \negthinspace }}\right] .  \notag
\end{gather}%
Again, it is clear from the structure of $t_{16}$ and $t_{16}^{T}$ that $%
s_{16}$ is a panmagic super-sudoku square as can be verified from the
defining formulas as before. Note that the elements of $s_{16}$ follow the
same permutation scheme as those of $s_{9}$ in (\ref{s9}).

A panmagic square can be constructed from $S_{16}\left( i,j\right) $ as%
\begin{gather}
M_{16}\left( i,j\right) =16S_{16}\left( i,j\right) +S_{16}\left(
i,17-j\right) \equiv m_{16},  \label{M16} \\
m_{16}=\left[ \!%
\begin{array}{cccccccccccccccc}
\text{{\scriptsize \negthinspace 15}} & \text{{\scriptsize \negthinspace 30}}
& \text{{\scriptsize \negthinspace 45}} & \text{{\scriptsize \negthinspace 60%
}} & \text{{\scriptsize \negthinspace 75}} & \text{{\scriptsize %
\negthinspace 90}} & \text{{\scriptsize \negthinspace 105}} & \text{%
{\scriptsize \negthinspace 120}} & \text{{\scriptsize \negthinspace 135}} & 
\text{{\scriptsize \negthinspace 150}} & \text{{\scriptsize \negthinspace 165%
}} & \text{{\scriptsize \negthinspace 180}} & \text{{\scriptsize %
\negthinspace 195}} & \text{{\scriptsize \negthinspace 210}} & \text{%
{\scriptsize \negthinspace 225}} & \text{{\scriptsize \negthinspace 240}} \\ 
\text{{\scriptsize \negthinspace 67}} & \text{{\scriptsize \negthinspace 82}}
& \text{{\scriptsize \negthinspace 97}} & \text{{\scriptsize \negthinspace
112}} & \text{{\scriptsize \negthinspace 143}} & \text{{\scriptsize %
\negthinspace 158}} & \text{{\scriptsize \negthinspace 173}} & \text{%
{\scriptsize \negthinspace 188}} & \text{{\scriptsize \negthinspace 203}} & 
\text{{\scriptsize \negthinspace 218}} & \text{{\scriptsize \negthinspace 233%
}} & \text{{\scriptsize \negthinspace 248}} & \text{{\scriptsize %
\negthinspace 7}} & \text{{\scriptsize \negthinspace 22}} & \text{%
{\scriptsize \negthinspace 37}} & \text{{\scriptsize \negthinspace 52}} \\ 
\text{{\scriptsize \negthinspace 135}} & \text{{\scriptsize \negthinspace 150%
}} & \text{{\scriptsize \negthinspace 165}} & \text{{\scriptsize %
\negthinspace 180}} & \text{{\scriptsize \negthinspace 195}} & \text{%
{\scriptsize \negthinspace 210}} & \text{{\scriptsize \negthinspace 225}} & 
\text{{\scriptsize \negthinspace 240}} & \text{{\scriptsize \negthinspace 15}%
} & \text{{\scriptsize \negthinspace 30}} & \text{{\scriptsize \negthinspace
45}} & \text{{\scriptsize \negthinspace 60}} & \text{{\scriptsize %
\negthinspace 75}} & \text{{\scriptsize \negthinspace 90}} & \text{%
{\scriptsize \negthinspace 105}} & \text{{\scriptsize \negthinspace 120}} \\ 
\text{{\scriptsize \negthinspace 203}} & \text{{\scriptsize \negthinspace 218%
}} & \text{{\scriptsize \negthinspace 233}} & \text{{\scriptsize %
\negthinspace 248}} & \text{{\scriptsize \negthinspace 7}} & \text{%
{\scriptsize \negthinspace 22}} & \text{{\scriptsize \negthinspace 37}} & 
\text{{\scriptsize \negthinspace 52}} & \text{{\scriptsize \negthinspace 67}}
& \text{{\scriptsize \negthinspace 82}} & \text{{\scriptsize \negthinspace 97%
}} & \text{{\scriptsize \negthinspace 112}} & \text{{\scriptsize %
\negthinspace 143}} & \text{{\scriptsize \negthinspace 158}} & \text{%
{\scriptsize \negthinspace 173}} & \text{{\scriptsize \negthinspace 188}} \\ 
\text{{\scriptsize \negthinspace 28}} & \text{{\scriptsize \negthinspace 47}}
& \text{{\scriptsize \negthinspace 62}} & \text{{\scriptsize \negthinspace 13%
}} & \text{{\scriptsize \negthinspace 88}} & \text{{\scriptsize %
\negthinspace 107}} & \text{{\scriptsize \negthinspace 122}} & \text{%
{\scriptsize \negthinspace 73}} & \text{{\scriptsize \negthinspace 148}} & 
\text{{\scriptsize \negthinspace 167}} & \text{{\scriptsize \negthinspace 182%
}} & \text{{\scriptsize \negthinspace 133}} & \text{{\scriptsize %
\negthinspace 208}} & \text{{\scriptsize \negthinspace 227}} & \text{%
{\scriptsize \negthinspace 242}} & \text{{\scriptsize \negthinspace 193}} \\ 
\text{{\scriptsize \negthinspace 80}} & \text{{\scriptsize \negthinspace 99}}
& \text{{\scriptsize \negthinspace 114}} & \text{{\scriptsize \negthinspace
65}} & \text{{\scriptsize \negthinspace 156}} & \text{{\scriptsize %
\negthinspace 175}} & \text{{\scriptsize \negthinspace 190}} & \text{%
{\scriptsize \negthinspace 141}} & \text{{\scriptsize \negthinspace 216}} & 
\text{{\scriptsize \negthinspace 235}} & \text{{\scriptsize \negthinspace 250%
}} & \text{{\scriptsize \negthinspace 201}} & \text{{\scriptsize %
\negthinspace 20}} & \text{{\scriptsize \negthinspace 39}} & \text{%
{\scriptsize \negthinspace 54}} & \text{{\scriptsize \negthinspace 5}} \\ 
\text{{\scriptsize \negthinspace 148}} & \text{{\scriptsize \negthinspace 167%
}} & \text{{\scriptsize \negthinspace 182}} & \text{{\scriptsize %
\negthinspace 133}} & \text{{\scriptsize \negthinspace 208}} & \text{%
{\scriptsize \negthinspace 227}} & \text{{\scriptsize \negthinspace 242}} & 
\text{{\scriptsize \negthinspace 193}} & \text{{\scriptsize \negthinspace 28}%
} & \text{{\scriptsize \negthinspace 47}} & \text{{\scriptsize \negthinspace
62}} & \text{{\scriptsize \negthinspace 13}} & \text{{\scriptsize %
\negthinspace 88}} & \text{{\scriptsize \negthinspace 107}} & \text{%
{\scriptsize \negthinspace 122}} & \text{{\scriptsize \negthinspace 73}} \\ 
\text{{\scriptsize \negthinspace 216}} & \text{{\scriptsize \negthinspace 235%
}} & \text{{\scriptsize \negthinspace 250}} & \text{{\scriptsize %
\negthinspace 201}} & \text{{\scriptsize \negthinspace 20}} & \text{%
{\scriptsize \negthinspace 39}} & \text{{\scriptsize \negthinspace 54}} & 
\text{{\scriptsize \negthinspace 5}} & \text{{\scriptsize \negthinspace 80}}
& \text{{\scriptsize \negthinspace 99}} & \text{{\scriptsize \negthinspace
114}} & \text{{\scriptsize \negthinspace 65}} & \text{{\scriptsize %
\negthinspace 156}} & \text{{\scriptsize \negthinspace 175}} & \text{%
{\scriptsize \negthinspace 190}} & \text{{\scriptsize \negthinspace 141}} \\ 
\text{{\scriptsize \negthinspace 45}} & \text{{\scriptsize \negthinspace 60}}
& \text{{\scriptsize \negthinspace 15}} & \text{{\scriptsize \negthinspace 30%
}} & \text{{\scriptsize \negthinspace 105}} & \text{{\scriptsize %
\negthinspace 120}} & \text{{\scriptsize \negthinspace 75}} & \text{%
{\scriptsize \negthinspace 90}} & \text{{\scriptsize \negthinspace 165}} & 
\text{{\scriptsize \negthinspace 180}} & \text{{\scriptsize \negthinspace 135%
}} & \text{{\scriptsize \negthinspace 150}} & \text{{\scriptsize %
\negthinspace 225}} & \text{{\scriptsize \negthinspace 240}} & \text{%
{\scriptsize \negthinspace 195}} & \text{{\scriptsize \negthinspace 210}} \\ 
\text{{\scriptsize \negthinspace 97}} & \text{{\scriptsize \negthinspace 112}%
} & \text{{\scriptsize \negthinspace 67}} & \text{{\scriptsize \negthinspace
82}} & \text{{\scriptsize \negthinspace 173}} & \text{{\scriptsize %
\negthinspace 188}} & \text{{\scriptsize \negthinspace 143}} & \text{%
{\scriptsize \negthinspace 158}} & \text{{\scriptsize \negthinspace 233}} & 
\text{{\scriptsize \negthinspace 248}} & \text{{\scriptsize \negthinspace 203%
}} & \text{{\scriptsize \negthinspace 218}} & \text{{\scriptsize %
\negthinspace 37}} & \text{{\scriptsize \negthinspace 52}} & \text{%
{\scriptsize \negthinspace 7}} & \text{{\scriptsize \negthinspace 22}} \\ 
\text{{\scriptsize \negthinspace 165}} & \text{{\scriptsize \negthinspace 180%
}} & \text{{\scriptsize \negthinspace 135}} & \text{{\scriptsize %
\negthinspace 150}} & \text{{\scriptsize \negthinspace 225}} & \text{%
{\scriptsize \negthinspace 240}} & \text{{\scriptsize \negthinspace 195}} & 
\text{{\scriptsize \negthinspace 210}} & \text{{\scriptsize \negthinspace 45}%
} & \text{{\scriptsize \negthinspace 60}} & \text{{\scriptsize \negthinspace
15}} & \text{{\scriptsize \negthinspace 30}} & \text{{\scriptsize %
\negthinspace 105}} & \text{{\scriptsize \negthinspace 120}} & \text{%
{\scriptsize \negthinspace 75}} & \text{{\scriptsize \negthinspace 90}} \\ 
\text{{\scriptsize \negthinspace 233}} & \text{{\scriptsize \negthinspace 248%
}} & \text{{\scriptsize \negthinspace 203}} & \text{{\scriptsize %
\negthinspace 218}} & \text{{\scriptsize \negthinspace 37}} & \text{%
{\scriptsize \negthinspace 52}} & \text{{\scriptsize \negthinspace 7}} & 
\text{{\scriptsize \negthinspace 22}} & \text{{\scriptsize \negthinspace 97}}
& \text{{\scriptsize \negthinspace 112}} & \text{{\scriptsize \negthinspace
67}} & \text{{\scriptsize \negthinspace 82}} & \text{{\scriptsize %
\negthinspace 173}} & \text{{\scriptsize \negthinspace 188}} & \text{%
{\scriptsize \negthinspace 143}} & \text{{\scriptsize \negthinspace 158}} \\ 
\text{{\scriptsize \negthinspace 62}} & \text{{\scriptsize \negthinspace 13}}
& \text{{\scriptsize \negthinspace 28}} & \text{{\scriptsize \negthinspace 47%
}} & \text{{\scriptsize \negthinspace 122}} & \text{{\scriptsize %
\negthinspace 73}} & \text{{\scriptsize \negthinspace 88}} & \text{%
{\scriptsize \negthinspace 107}} & \text{{\scriptsize \negthinspace 182}} & 
\text{{\scriptsize \negthinspace 133}} & \text{{\scriptsize \negthinspace 148%
}} & \text{{\scriptsize \negthinspace 167}} & \text{{\scriptsize %
\negthinspace 242}} & \text{{\scriptsize \negthinspace 193}} & \text{%
{\scriptsize \negthinspace 208}} & \text{{\scriptsize \negthinspace 227}} \\ 
\text{{\scriptsize \negthinspace 114}} & \text{{\scriptsize \negthinspace 65}%
} & \text{{\scriptsize \negthinspace 80}} & \text{{\scriptsize \negthinspace
99}} & \text{{\scriptsize \negthinspace 190}} & \text{{\scriptsize %
\negthinspace 141}} & \text{{\scriptsize \negthinspace 156}} & \text{%
{\scriptsize \negthinspace 175}} & \text{{\scriptsize \negthinspace 250}} & 
\text{{\scriptsize \negthinspace 201}} & \text{{\scriptsize \negthinspace 216%
}} & \text{{\scriptsize \negthinspace 235}} & \text{{\scriptsize %
\negthinspace 54}} & \text{{\scriptsize \negthinspace 5}} & \text{%
{\scriptsize \negthinspace 20}} & \text{{\scriptsize \negthinspace 39}} \\ 
\text{{\scriptsize \negthinspace 182}} & \text{{\scriptsize \negthinspace 133%
}} & \text{{\scriptsize \negthinspace 148}} & \text{{\scriptsize %
\negthinspace 167}} & \text{{\scriptsize \negthinspace 242}} & \text{%
{\scriptsize \negthinspace 193}} & \text{{\scriptsize \negthinspace 208}} & 
\text{{\scriptsize \negthinspace 227}} & \text{{\scriptsize \negthinspace 62}%
} & \text{{\scriptsize \negthinspace 13}} & \text{{\scriptsize \negthinspace
28}} & \text{{\scriptsize \negthinspace 47}} & \text{{\scriptsize %
\negthinspace 122}} & \text{{\scriptsize \negthinspace 73}} & \text{%
{\scriptsize \negthinspace 88}} & \text{{\scriptsize \negthinspace 107}} \\ 
\text{{\scriptsize \negthinspace 250}} & \text{{\scriptsize \negthinspace 201%
}} & \text{{\scriptsize \negthinspace 216}} & \text{{\scriptsize %
\negthinspace 235}} & \text{{\scriptsize \negthinspace 54}} & \text{%
{\scriptsize \negthinspace 5}} & \text{{\scriptsize \negthinspace 20}} & 
\text{{\scriptsize \negthinspace 39}} & \text{{\scriptsize \negthinspace 114}%
} & \text{{\scriptsize \negthinspace 65}} & \text{{\scriptsize \negthinspace
80}} & \text{{\scriptsize \negthinspace 99}} & \text{{\scriptsize %
\negthinspace 190}} & \text{{\scriptsize \negthinspace 141}} & \text{%
{\scriptsize \negthinspace 156}} & \text{{\scriptsize \negthinspace 175}}%
\end{array}%
\!\!\right] .  \notag
\end{gather}%
However, this square is not natural as it contains many duplicate elements.
As before, it can be verified that that $m_{16}$ is a panmagic square whose
4 by 4 subsquares all add to $m=2040\,.$ Higher order-$k^{2}$ pandiagonal
super-sudoku squares can be generated and studied in a similar manner. For
such squares, we conjecture that a panmagic square formed from (\ref{Mn}) is
natural only when $k$ is odd.$\smallskip $

\begin{center}
{\large Knut Vik Sudoku Squares}
\end{center}

\noindent According to (\ref{Kij}), an order-25 Knut Vik sudoku square (the
lowest possible order) is generated by%
\begin{gather}
V_{25}(i,j)=5\left( \left( 1-j+2\left( i-1\right) \right) \func{mod}5\right) 
\notag \\
+\left( \left( -2(j-1)+2\left( i-1\right) +2\left\lfloor \left( i-1\right)
/5\right\rfloor +\left\lfloor \left( j-1\right) /5\right\rfloor \right) 
\func{mod}5\right) ,  \label{V25} \\
i,j=1,2,\ldots ,25\,.  \notag
\end{gather}%
We also can express $V_{25}$ in terms of order-$5$ submatrices $A_{k}$ as 
\begin{equation}
V_{25}=\left[ 
\begin{array}{rrrrr}
Z_{0} & Z_{1} & Z_{2} & Z_{3} & Z_{4} \\ 
Z_{2} & Z_{3} & Z_{4} & Z_{0} & Z_{1} \\ 
Z_{4} & Z_{0} & Z_{1} & Z_{2} & Z_{3} \\ 
Z_{1} & Z_{2} & Z_{3} & Z_{4} & Z_{0} \\ 
Z_{3} & Z_{4} & Z_{0} & Z_{1} & Z_{2}%
\end{array}%
\right] ,  \label{VZ25}
\end{equation}%
where, as indicated by (\ref{MauxG}), (\ref{Mt5}), and (\ref{M5sub}), 
\begin{eqnarray}
Z_{i} &=&5\left[ 
\begin{array}{rrrrr}
0 & 4 & 3 & 2 & 1 \\ 
2 & 1 & 0 & 4 & 3 \\ 
4 & 3 & 2 & 1 & 0 \\ 
1 & 0 & 4 & 3 & 2 \\ 
3 & 2 & 1 & 0 & 4%
\end{array}%
\right] +\left[ 
\begin{array}{ccccc}
i & i+3 & i+1 & i+4 & i+2 \\ 
i+2 & i & i+3 & i+1 & i+4 \\ 
i+4 & i+2 & i & i+3 & i+1 \\ 
i+1 & i+4 & i+2 & i & i+3 \\ 
i+3 & i+1 & i+4 & i+2 & i%
\end{array}%
\right] \func{mod}5, \\
i &=&0,1,\ldots ,4\,.  \notag
\end{eqnarray}%
The structure of these submatrices shows why $V_{25}$ is a Knut Vik sudoku
square. These properties also can be verified by evaluating (\ref{Lat}), (%
\ref{KV}), and (\ref{Sud}). However, evaluation of (\ref{KSK}) shows that $%
V_{25}$ is not a super-sudoku square. By (\ref{V25}) or (\ref{VZ25}) we have
the following order-25 Knut Vik sudoku square:\smallskip 

\noindent $%
\begin{tabular}[t]{|p{0.04in}p{0.04in}p{0.04in}p{0.04in}p{0.04in}|p{0.04in}p{0.04in}p{0.04in}p{0.04in}p{0.04in}|p{0.04in}p{0.04in}p{0.04in}p{0.04in}p{0.04in}|p{0.04in}p{0.04in}p{0.04in}p{0.04in}p{0.04in}|p{0.04in}p{0.04in}p{0.04in}p{0.04in}p{0.04in}|}
\hline
{\scriptsize 0} & {\scriptsize 23} & {\scriptsize 16} & {\scriptsize 14} & 
{\scriptsize 7} & {\scriptsize 1} & {\scriptsize 24} & {\scriptsize 17} & 
{\scriptsize 10} & {\scriptsize 8} & {\scriptsize 2} & {\scriptsize 20} & 
{\scriptsize 18} & {\scriptsize 11} & {\scriptsize 9} & {\scriptsize 3} & 
{\scriptsize 21} & {\scriptsize 19} & {\scriptsize 12} & {\scriptsize 5} & 
{\scriptsize 4} & {\scriptsize 22} & {\scriptsize 15} & {\scriptsize 13} & 
{\scriptsize 6} \\ 
{\scriptsize 12} & {\scriptsize 5} & {\scriptsize 3} & {\scriptsize 21} & 
{\scriptsize 19} & {\scriptsize 13} & {\scriptsize 6} & {\scriptsize 4} & 
{\scriptsize 22} & {\scriptsize 15} & {\scriptsize 14} & {\scriptsize 7} & 
{\scriptsize 0} & {\scriptsize 23} & {\scriptsize 16} & {\scriptsize 10} & 
{\scriptsize 8} & {\scriptsize 1} & {\scriptsize 24} & {\scriptsize 17} & 
{\scriptsize 11} & {\scriptsize 9} & {\scriptsize 2} & {\scriptsize 20} & 
{\scriptsize 18} \\ 
{\scriptsize 24} & {\scriptsize 17} & {\scriptsize 10} & {\scriptsize 8} & 
{\scriptsize 1} & {\scriptsize 20} & {\scriptsize 18} & {\scriptsize 11} & 
{\scriptsize 9} & {\scriptsize 2} & {\scriptsize 21} & {\scriptsize 19} & 
{\scriptsize 12} & {\scriptsize 5} & {\scriptsize 3} & {\scriptsize 22} & 
{\scriptsize 15} & {\scriptsize 13} & {\scriptsize 6} & {\scriptsize 4} & 
{\scriptsize 23} & {\scriptsize 16} & {\scriptsize 14} & {\scriptsize 7} & 
{\scriptsize 0} \\ 
{\scriptsize 6} & {\scriptsize 4} & {\scriptsize 22} & {\scriptsize 15} & 
{\scriptsize 13} & {\scriptsize 7} & {\scriptsize 0} & {\scriptsize 23} & 
{\scriptsize 16} & {\scriptsize 14} & {\scriptsize 8} & {\scriptsize 1} & 
{\scriptsize 24} & {\scriptsize 17} & {\scriptsize 10} & {\scriptsize 9} & 
{\scriptsize 2} & {\scriptsize 20} & {\scriptsize 18} & {\scriptsize 11} & 
{\scriptsize 5} & {\scriptsize 3} & {\scriptsize 21} & {\scriptsize 19} & 
{\scriptsize 12} \\ 
{\scriptsize 18} & {\scriptsize 11} & {\scriptsize 9} & {\scriptsize 2} & 
{\scriptsize 20} & {\scriptsize 19} & {\scriptsize 12} & {\scriptsize 5} & 
{\scriptsize 3} & {\scriptsize 21} & {\scriptsize 15} & {\scriptsize 13} & 
{\scriptsize 6} & {\scriptsize 4} & {\scriptsize 22} & {\scriptsize 16} & 
{\scriptsize 14} & {\scriptsize 7} & {\scriptsize 0} & {\scriptsize 23} & 
{\scriptsize 17} & {\scriptsize 10} & {\scriptsize 8} & {\scriptsize 1} & 
{\scriptsize 24} \\ \hline
{\scriptsize 2} & {\scriptsize 20} & {\scriptsize 18} & {\scriptsize 11} & 
{\scriptsize 9} & {\scriptsize 3} & {\scriptsize 21} & {\scriptsize 19} & 
{\scriptsize 12} & {\scriptsize 5} & {\scriptsize 4} & {\scriptsize 22} & 
{\scriptsize 15} & {\scriptsize 13} & {\scriptsize 6} & {\scriptsize 0} & 
{\scriptsize 23} & {\scriptsize 16} & {\scriptsize 14} & {\scriptsize 7} & 
{\scriptsize 1} & {\scriptsize 24} & {\scriptsize 17} & {\scriptsize 10} & 
{\scriptsize 8} \\ 
{\scriptsize 14} & {\scriptsize 7} & {\scriptsize 0} & {\scriptsize 23} & 
{\scriptsize 16} & {\scriptsize 10} & {\scriptsize 8} & {\scriptsize 1} & 
{\scriptsize 24} & {\scriptsize 17} & {\scriptsize 11} & {\scriptsize 9} & 
{\scriptsize 2} & {\scriptsize 20} & {\scriptsize 18} & {\scriptsize 12} & 
{\scriptsize 5} & {\scriptsize 3} & {\scriptsize 21} & {\scriptsize 19} & 
{\scriptsize 13} & {\scriptsize 6} & {\scriptsize 4} & {\scriptsize 22} & 
{\scriptsize 15} \\ 
{\scriptsize 21} & {\scriptsize 19} & {\scriptsize 12} & {\scriptsize 5} & 
{\scriptsize 3} & {\scriptsize 22} & {\scriptsize 15} & {\scriptsize 13} & 
{\scriptsize 6} & {\scriptsize 4} & {\scriptsize 23} & {\scriptsize 16} & 
{\scriptsize 14} & {\scriptsize 7} & {\scriptsize 0} & {\scriptsize 24} & 
{\scriptsize 17} & {\scriptsize 10} & {\scriptsize 8} & {\scriptsize 1} & 
{\scriptsize 20} & {\scriptsize 18} & {\scriptsize 11} & {\scriptsize 9} & 
{\scriptsize 2} \\ 
{\scriptsize 8} & {\scriptsize 1} & {\scriptsize 24} & {\scriptsize 17} & 
{\scriptsize 10} & {\scriptsize 9} & {\scriptsize 2} & {\scriptsize 20} & 
{\scriptsize 18} & {\scriptsize 11} & {\scriptsize 5} & {\scriptsize 3} & 
{\scriptsize 21} & {\scriptsize 19} & {\scriptsize 12} & {\scriptsize 6} & 
{\scriptsize 4} & {\scriptsize 22} & {\scriptsize 15} & {\scriptsize 13} & 
{\scriptsize 7} & {\scriptsize 0} & {\scriptsize 23} & {\scriptsize 16} & 
{\scriptsize 14} \\ 
{\scriptsize 15} & {\scriptsize 13} & {\scriptsize 6} & {\scriptsize 4} & 
{\scriptsize 22} & {\scriptsize 16} & {\scriptsize 14} & {\scriptsize 7} & 
{\scriptsize 0} & {\scriptsize 23} & {\scriptsize 17} & {\scriptsize 10} & 
{\scriptsize 8} & {\scriptsize 1} & {\scriptsize 24} & {\scriptsize 18} & 
{\scriptsize 11} & {\scriptsize 9} & {\scriptsize 2} & {\scriptsize 20} & 
{\scriptsize 19} & {\scriptsize 12} & {\scriptsize 5} & {\scriptsize 3} & 
{\scriptsize 21} \\ \hline
{\scriptsize 4} & {\scriptsize 22} & {\scriptsize 15} & {\scriptsize 13} & 
{\scriptsize 6} & {\scriptsize 0} & {\scriptsize 23} & {\scriptsize 16} & 
{\scriptsize 14} & {\scriptsize 7} & {\scriptsize 1} & {\scriptsize 24} & 
{\scriptsize 17} & {\scriptsize 10} & {\scriptsize 8} & {\scriptsize 2} & 
{\scriptsize 20} & {\scriptsize 18} & {\scriptsize 11} & {\scriptsize 9} & 
{\scriptsize 3} & {\scriptsize 21} & {\scriptsize 19} & {\scriptsize 12} & 
{\scriptsize 5} \\ 
{\scriptsize 11} & {\scriptsize 9} & {\scriptsize 2} & {\scriptsize 20} & 
{\scriptsize 18} & {\scriptsize 12} & {\scriptsize 5} & {\scriptsize 3} & 
{\scriptsize 21} & {\scriptsize 19} & {\scriptsize 13} & {\scriptsize 6} & 
{\scriptsize 4} & {\scriptsize 22} & {\scriptsize 15} & {\scriptsize 14} & 
{\scriptsize 7} & {\scriptsize 0} & {\scriptsize 23} & {\scriptsize 16} & 
{\scriptsize 10} & {\scriptsize 8} & {\scriptsize 1} & {\scriptsize 24} & 
{\scriptsize 17} \\ 
{\scriptsize 23} & {\scriptsize 16} & {\scriptsize 14} & {\scriptsize 7} & 
{\scriptsize 0} & {\scriptsize 24} & {\scriptsize 17} & {\scriptsize 10} & 
{\scriptsize 8} & {\scriptsize 1} & {\scriptsize 20} & {\scriptsize 18} & 
{\scriptsize 11} & {\scriptsize 9} & {\scriptsize 2} & {\scriptsize 21} & 
{\scriptsize 19} & {\scriptsize 12} & {\scriptsize 5} & {\scriptsize 3} & 
{\scriptsize 22} & {\scriptsize 15} & {\scriptsize 13} & {\scriptsize 6} & 
{\scriptsize 4} \\ 
{\scriptsize 5} & {\scriptsize 3} & {\scriptsize 21} & {\scriptsize 19} & 
{\scriptsize 12} & {\scriptsize 6} & {\scriptsize 4} & {\scriptsize 22} & 
{\scriptsize 15} & {\scriptsize 13} & {\scriptsize 7} & {\scriptsize 0} & 
{\scriptsize 23} & {\scriptsize 16} & {\scriptsize 14} & {\scriptsize 8} & 
{\scriptsize 1} & {\scriptsize 24} & {\scriptsize 17} & {\scriptsize 10} & 
{\scriptsize 9} & {\scriptsize 2} & {\scriptsize 20} & {\scriptsize 18} & 
{\scriptsize 11} \\ 
{\scriptsize 17} & {\scriptsize 10} & {\scriptsize 8} & {\scriptsize 1} & 
{\scriptsize 24} & {\scriptsize 18} & {\scriptsize 11} & {\scriptsize 9} & 
{\scriptsize 2} & {\scriptsize 20} & {\scriptsize 19} & {\scriptsize 12} & 
{\scriptsize 5} & {\scriptsize 3} & {\scriptsize 21} & {\scriptsize 15} & 
{\scriptsize 13} & {\scriptsize 6} & {\scriptsize 4} & {\scriptsize 22} & 
{\scriptsize 16} & {\scriptsize 14} & {\scriptsize 7} & {\scriptsize 0} & 
{\scriptsize 23} \\ \hline
{\scriptsize 1} & {\scriptsize 24} & {\scriptsize 17} & {\scriptsize 10} & 
{\scriptsize 8} & {\scriptsize 2} & {\scriptsize 20} & {\scriptsize 18} & 
{\scriptsize 11} & {\scriptsize 9} & {\scriptsize 3} & {\scriptsize 21} & 
{\scriptsize 19} & {\scriptsize 12} & {\scriptsize 5} & {\scriptsize 4} & 
{\scriptsize 22} & {\scriptsize 15} & {\scriptsize 13} & {\scriptsize 6} & 
{\scriptsize 0} & {\scriptsize 23} & {\scriptsize 16} & {\scriptsize 14} & 
{\scriptsize 7} \\ 
{\scriptsize 13} & {\scriptsize 6} & {\scriptsize 4} & {\scriptsize 22} & 
{\scriptsize 15} & {\scriptsize 14} & {\scriptsize 7} & {\scriptsize 0} & 
{\scriptsize 23} & {\scriptsize 16} & {\scriptsize 10} & {\scriptsize 8} & 
{\scriptsize 1} & {\scriptsize 24} & {\scriptsize 17} & {\scriptsize 11} & 
{\scriptsize 9} & {\scriptsize 2} & {\scriptsize 20} & {\scriptsize 18} & 
{\scriptsize 12} & {\scriptsize 5} & {\scriptsize 3} & {\scriptsize 21} & 
{\scriptsize 19} \\ 
{\scriptsize 20} & {\scriptsize 18} & {\scriptsize 11} & {\scriptsize 9} & 
{\scriptsize 2} & {\scriptsize 21} & {\scriptsize 19} & {\scriptsize 12} & 
{\scriptsize 5} & {\scriptsize 3} & {\scriptsize 22} & {\scriptsize 15} & 
{\scriptsize 13} & {\scriptsize 6} & {\scriptsize 4} & {\scriptsize 23} & 
{\scriptsize 16} & {\scriptsize 14} & {\scriptsize 7} & {\scriptsize 0} & 
{\scriptsize 24} & {\scriptsize 17} & {\scriptsize 10} & {\scriptsize 8} & 
{\scriptsize 1} \\ 
{\scriptsize 7} & {\scriptsize 0} & {\scriptsize 23} & {\scriptsize 16} & 
{\scriptsize 14} & {\scriptsize 8} & {\scriptsize 1} & {\scriptsize 24} & 
{\scriptsize 17} & {\scriptsize 10} & {\scriptsize 9} & {\scriptsize 2} & 
{\scriptsize 20} & {\scriptsize 18} & {\scriptsize 11} & {\scriptsize 5} & 
{\scriptsize 3} & {\scriptsize 21} & {\scriptsize 19} & {\scriptsize 12} & 
{\scriptsize 6} & {\scriptsize 4} & {\scriptsize 22} & {\scriptsize 15} & 
{\scriptsize 13} \\ 
{\scriptsize 19} & {\scriptsize 12} & {\scriptsize 5} & {\scriptsize 3} & 
{\scriptsize 21} & {\scriptsize 15} & {\scriptsize 13} & {\scriptsize 6} & 
{\scriptsize 4} & {\scriptsize 22} & {\scriptsize 16} & {\scriptsize 14} & 
{\scriptsize 7} & {\scriptsize 0} & {\scriptsize 23} & {\scriptsize 17} & 
{\scriptsize 10} & {\scriptsize 8} & {\scriptsize 1} & {\scriptsize 24} & 
{\scriptsize 18} & {\scriptsize 11} & {\scriptsize 9} & {\scriptsize 2} & 
{\scriptsize 20} \\ \hline
{\scriptsize 3} & {\scriptsize 21} & {\scriptsize 19} & {\scriptsize 12} & 
{\scriptsize 5} & {\scriptsize 4} & {\scriptsize 22} & {\scriptsize 15} & 
{\scriptsize 13} & {\scriptsize 6} & {\scriptsize 0} & {\scriptsize 23} & 
{\scriptsize 16} & {\scriptsize 14} & {\scriptsize 7} & {\scriptsize 1} & 
{\scriptsize 24} & {\scriptsize 17} & {\scriptsize 10} & {\scriptsize 8} & 
{\scriptsize 2} & {\scriptsize 20} & {\scriptsize 18} & {\scriptsize 11} & 
{\scriptsize 9} \\ 
{\scriptsize 10} & {\scriptsize 8} & {\scriptsize 1} & {\scriptsize 24} & 
{\scriptsize 17} & {\scriptsize 11} & {\scriptsize 9} & {\scriptsize 2} & 
{\scriptsize 20} & {\scriptsize 18} & {\scriptsize 12} & {\scriptsize 5} & 
{\scriptsize 3} & {\scriptsize 21} & {\scriptsize 19} & {\scriptsize 13} & 
{\scriptsize 6} & {\scriptsize 4} & {\scriptsize 22} & {\scriptsize 15} & 
{\scriptsize 14} & {\scriptsize 7} & {\scriptsize 0} & {\scriptsize 23} & 
{\scriptsize 16} \\ 
{\scriptsize 22} & {\scriptsize 15} & {\scriptsize 13} & {\scriptsize 6} & 
{\scriptsize 4} & {\scriptsize 23} & {\scriptsize 16} & {\scriptsize 14} & 
{\scriptsize 7} & {\scriptsize 0} & {\scriptsize 24} & {\scriptsize 17} & 
{\scriptsize 10} & {\scriptsize 8} & {\scriptsize 1} & {\scriptsize 20} & 
{\scriptsize 18} & {\scriptsize 11} & {\scriptsize 9} & {\scriptsize 2} & 
{\scriptsize 21} & {\scriptsize 19} & {\scriptsize 12} & {\scriptsize 5} & 
{\scriptsize 3} \\ 
{\scriptsize 9} & {\scriptsize 2} & {\scriptsize 20} & {\scriptsize 18} & 
{\scriptsize 11} & {\scriptsize 5} & {\scriptsize 3} & {\scriptsize 21} & 
{\scriptsize 19} & {\scriptsize 12} & {\scriptsize 6} & {\scriptsize 4} & 
{\scriptsize 22} & {\scriptsize 15} & {\scriptsize 13} & {\scriptsize 7} & 
{\scriptsize 0} & {\scriptsize 23} & {\scriptsize 16} & {\scriptsize 14} & 
{\scriptsize 8} & {\scriptsize 1} & {\scriptsize 24} & {\scriptsize 17} & 
{\scriptsize 10} \\ 
{\scriptsize 16} & {\scriptsize 14} & {\scriptsize 7} & {\scriptsize 0} & 
{\scriptsize 23} & {\scriptsize 17} & {\scriptsize 10} & {\scriptsize 8} & 
{\scriptsize 1} & {\scriptsize 24} & {\scriptsize 18} & {\scriptsize 11} & 
{\scriptsize 9} & {\scriptsize 2} & {\scriptsize 20} & {\scriptsize 19} & 
{\scriptsize 12} & {\scriptsize 5} & {\scriptsize 3} & {\scriptsize 21} & 
{\scriptsize 15} & {\scriptsize 13} & {\scriptsize 6} & {\scriptsize 4} & 
{\scriptsize 22} \\ \hline
\end{tabular}%
\bigskip $

\noindent As noted after (\ref{Kij}), this square is essentially the same as
the one given by Boyer \cite{BOY}. It differs only by a shift of columns and
the addition of $1$ to each element.

Knut Vik sudoku squares $V_{n}$ of higher order $n=k^{2}$ ($k$ not divisible
by 2 or 3) can be constructed in a similar manner as for $n=25$. An order-49
of this type has the form\footnote{%
The full 49 by 49 square can be constructed from (\ref{Kij}) in Excel%
\copyright .}%
\begin{equation}
V_{49}=\left[ 
\begin{array}{ccccccc}
W_{0} & W_{1} & W_{2} & W_{3} & W_{4} & W_{5} & W_{6} \\ 
W_{2} & W_{3} & W_{4} & W_{5} & W_{6} & W_{0} & W_{1} \\ 
W_{4} & W_{5} & W_{6} & W_{0} & W_{1} & W_{2} & W_{3} \\ 
W_{6} & W_{0} & W_{1} & W_{2} & W_{3} & W_{4} & W_{5} \\ 
W_{1} & W_{2} & W_{3} & W_{4} & W_{5} & W_{6} & W_{0} \\ 
W_{3} & W_{4} & W_{5} & W_{6} & W_{0} & W_{1} & W_{2} \\ 
W_{5} & W_{6} & W_{0} & W_{1} & W_{2} & W_{3} & W_{4}%
\end{array}%
\right] ,
\end{equation}%
where the submatrices $W_{i}$ are given by (\ref{Kij}) as%
\begin{eqnarray}
W_{i} &=&7\left[ 
\begin{array}{ccccccc}
0 & 6 & 5 & 4 & 3 & 2 & 1 \\ 
2 & 1 & 0 & 6 & 5 & 4 & 3 \\ 
4 & 3 & 2 & 1 & 0 & 6 & 5 \\ 
6 & 5 & 4 & 3 & 2 & 1 & 0 \\ 
1 & 0 & 6 & 5 & 4 & 3 & 2 \\ 
3 & 2 & 1 & 0 & 6 & 5 & 4 \\ 
5 & 4 & 3 & 2 & 1 & 0 & 6%
\end{array}%
\right]   \notag \\
&&+\left[ 
\begin{array}{ccccccc}
i & i+5 & i+3 & i+1 & i+6 & i+4 & i+2 \\ 
i+2 & i & i+5 & i+3 & i+1 & i+6 & i+4 \\ 
i+4 & i+2 & i & i+5 & i+3 & i+1 & i+6 \\ 
i+6 & i+4 & i+2 & i & i+5 & i+3 & i+1 \\ 
i+1 & i+6 & i+4 & i+2 & i & i+5 & i+3 \\ 
i+3 & i+1 & i+6 & i+4 & i+2 & i & i+5 \\ 
i+5 & i+3 & i+1 & i+6 & i+4 & i+2 & i%
\end{array}%
\right] \func{mod}7, \\
i &=&0,1,\ldots ,6,  \notag
\end{eqnarray}%
i.e.,%
\begin{eqnarray}
W_{0} &=&\left[ 
\begin{array}{ccccccc}
0 & 47 & 38 & 29 & 27 & 18 & 9 \\ 
16 & 7 & 5 & 45 & 36 & 34 & 25 \\ 
32 & 23 & 14 & 12 & 3 & 43 & 41 \\ 
48 & 39 & 30 & 21 & 19 & 10 & 1 \\ 
8 & 6 & 46 & 37 & 28 & 26 & 17 \\ 
24 & 15 & 13 & 4 & 44 & 35 & 33 \\ 
40 & 31 & 22 & 20 & 11 & 2 & 42%
\end{array}%
\allowbreak \right] ,\cdots ,  \notag \\
W_{6} &=&\left[ 
\begin{array}{ccccccc}
6 & 46 & 37 & 28 & 26 & 17 & 8 \\ 
15 & 13 & 4 & 44 & 35 & 33 & 24 \\ 
31 & 22 & 20 & 11 & 2 & 42 & 40 \\ 
47 & 38 & 29 & 27 & 18 & 9 & 0 \\ 
7 & 5 & 45 & 36 & 34 & 25 & 16 \\ 
23 & 14 & 12 & 3 & 43 & 41 & 32 \\ 
39 & 30 & 21 & 19 & 10 & 1 & 48%
\end{array}%
\right] .
\end{eqnarray}

The algorithm (\ref{Kij}) for Knut Vik sudoku squares can be compared with
the following algorithm for Knut Vik squares of prime order $n>3\ $due to
Euler as noted by Hedayat \cite{HED}:%
\begin{equation}
\hat{V}_{n}\left( i,j\right) =\left( \lambda i+j\right) \func{mod}n,\quad
\lambda \neq 0,\ 1,\ n-1\,.  \label{Vhat}
\end{equation}%
For nonprime $n$ not divisible by 2 or 3, Hedayat gives additional
restrictions on $\lambda $. The squares generated by (\ref{Vhat}) are
cyclic, whereas $V_{n}$ generated by (\ref{Kij}) are not. 

The $k$ by $k$ sudoku subsquares of $V_{n}$ are natural and semi-magic. In
addition, a natural panmagic square follows from $V_{n}$ using an equation
of the form (\ref{Mn}) since $V_{n}$ and $V_{n}R$ are orthogonal as can be
shown by an extension of Br\'{e}e's orthogonality criterion as derived by
Nordgren \cite{NOR1}. For this extension, note that in each column of sudoku
squares each number follows a path that starts in the first column and
progresses down $k+1$, right two (with wraparound in the column).
Furthermore, the sudoku squares themselves follow a knight's move pattern.
Therefore, as in \cite{NOR1}, analysis shows that the paths for each number
pair in $V_{n}$ and $V_{n}R$ intersect only once and orthogonality follows.

\end{document}